\input amstex
\documentstyle{amsppt}
\topmatter
\title Cohomology of complete intersections in toric varieties
\endtitle
\author Anvar R. Mavlyutov\endauthor
\affil University of Massachusetts at Amherst\endaffil
\address Department of Mathematics and Statistics,
Univ. of Massachusetts, Amherst, MA, 01003, USA
\endaddress
\email anvar\@math.umass.edu\endemail

\abstract We explicitly describe cohomology of complete intersections
in compact simplicial toric varieties.\endabstract
\endtopmatter

\document
In this  paper we will study intersections of hypersurfaces 
in compact simplicial toric varieties $\bold P_\Sigma$.
The main purpose is to relate naturally the Hodge structure
of a complete intersection $X_{f_1}\cap\ldots\cap X_{f_s}$ in 
$\bold P_\Sigma$ to a graded ring.
Originally this idea appears in  [Gr], [St], [Dol], [PS].
The case of a hypersurface in a toric variety has been treated
in [BC]. Also the Hodge structure of complete intersections
in a projective space was described in [Te], [Ko], [L], [Di], [Na]. 
The common approach was to reduce  studying of the Hodge structure 
on a complete intersection to  studying of the Hodge structure
on a hypersurface in a higher dimensional projective variety.
This is the idea of a ``Cayley trick".
About a Cayley trick in the toric context see [GKZ], [DK], [BB].
A special case of a complete intersection (when it is empty)
in a complete simplicial toric variety was elaborated in [CCD].
The basic references on toric varieties are [F1], [O], [Da], [C]. 

The paper is organized as follows:

Section 1 establishes notation and studies cohomology of subvarieties
in a complete simplicial toric variety. 
In Section 2 we describe a Cayley
trick for toric varieties.
In Section 3 we prove the main result where we relate 
the Hodge components $H^{d-s-p,p}(X_{f_1}\cap\ldots\cap X_{f_s})$ in the 
middle cohomology group to homogeneous components of a graded ring.
Section 4 treats
a special case of complete intersections: 
 a nondegenerate intersection.

{\it Acknowledgment.} I would like to thank D. Cox for his advice
and useful comments. 
    
{\bf 1. Quasi-smooth intersections.}
We first fix some notation.                                     
Let $M$ be a lattice of rank $d$, $N=\text{Hom}(M,Z)$ the 
dual lattice; $M_{\bold R}$ (resp. $N_{\bold R})$ denotes
the $\bold R$-scalar extension of $M$ (resp. of $N$).
Let $\Sigma$ be a rational simplicial complete  $d$-dimensional fan 
in $N_{\bold R}$ [BC], $\bold P_{\Sigma}$ a complete simplicial 
toric variety associated with this fan.  
                                  
Such a toric variety can be described as a geometric quotient [C].
Let
$S(\Sigma)=\bold C[x_1,\dots,x_n]$ be the polynomial ring over $\bold C$
with variables $x_1,\dots,x_n$ 
corresponding to the integral generators $e_1,\dots,e_n$
of the  1-dimensional cones  of
$\Sigma$.
For $\sigma\in\Sigma$ let $\hat{x}_{\sigma}=\prod_{e_i\notin\sigma}x_i$,
and let $B(\Sigma)=\langle\hat x_\sigma :\sigma\in\Sigma\rangle\subset S$ 
be the ideal
generated by the $\hat x_\sigma$'s. This ideal gives the variety
$Z(\Sigma)=\bold V(B(\Sigma))\subset \bold A^n$.
The toric variety $\bold P=\bold P_\Sigma$ will be a geometric quotient
of $U(\Sigma):=\bold A^n\setminus Z(\Sigma)$ 
by the group $\bold D:=\text{Hom}_{\bold Z}(A_{d-1}(\bold P), \bold C^*)$,
where $A_{d-1}(\bold P)$ is the Chow group of Weil divisors modulo 
rational equivalence.

Each variable $x_i$ in the coordinate ring $S(\Sigma)$
 corresponds to a torus-invariant
irreducible divisor $D_i$ of $\bold P$. As in [C], we grade $S=S(\Sigma)$ 
by assigning to a monomial $\prod^{n}_{i=1}x_{i}^{a_i}$ its degree 
$[\sum^{n}_{i=1}a_i D_i]\in A_{d-1}(\bold P)$.
A polynomial $f$ in the graded piece $S_\alpha$ corresponding to 
$\alpha\in A_{d-1}(\bold P)$
is said to be $\bold D$-homogeneous of degree $\alpha$.
                                                
Let $f_1,\dots,f_s$ be $\bold D$-homogeneous polynomials.
They define a zero set $\bold V(f_1,\dots,$ $f_s)\subset\bold A^n$,
 moreover 
$\bold V(f_1,\dots,f_s)\cap U(\Sigma)$ is stable under the action of $\bold D$
and hence descends to a closed subset $X\subset\bold P$,
because $\bold P$ is a geometric quotient.
                  
{\bf Definition 1.1.}                                   
We say that $X$ is a {\it quasi-smooth intersection} if
$\bold V(f_1,\dots,f_s)\cap U(\Sigma)$ is either empty or
 a smooth subvariety of
codimension $s$ in $U(\Sigma)$.
                             
{\bf Remark 1.2.} 
This notion generalizes a nonsingular complete intersection in a projective
space.
Notice that since the $(n-d)$-dimensional group $\bold D$ has only 
zero dimensional stabilizers [BC], 
$X$ is of pure dimension $d-s$ or empty.

We can now relate this notion to a V-submanifold (see Definition 3.2 in [BC]).

{\bf Proposition 1.3.}
{\it If $X\subset\bold P$ is a closed subset of codimension $s$
defined by 
$\bold D$-homogeneous polynomials 
$f_1,\dots,f_s$, then $X$ is a quasi-smooth intersection 
if and only if $X$ is a V-submanifold of $\bold P$.}

The proof of this is very similar to the proof of the Proposition 3.5 in [BC].

The next result is a Lefschetz-type theorem.
                 
{\bf Proposition 1.4.}
{\it                  
Let $X\subset\bold P$ be a closed subset,
defined by $\bold D$-homogeneous polynomials
$f_1,\dots,f_s$, in a complete simplicial toric variety $\bold P$.
If $f_1,\dots,f_s\in B(\Sigma)$, then the natural map 
$i^{*}: H^{i}(\bold P)\rightarrow H^{i}(X)$ is an 
isomorphism for $i<d-s$ and an injection for $i=d-s$.
In particular, 
this is valid if $X$ is an intersection of ample hypersurfaces.}
                                    
{\it Proof.}
We can present $X=X_{f_1}\cap\ldots\cap X_{f_s}$, where $X_{f_i}\subset\bold P$ 
is a hypersurface defined by $f_i$.
As it was shown in the proof of the 
Proposition 10.8 [BC], if $f\in B(\Sigma)$ then
$\bold P\setminus X_{f}=(\bold A^n\setminus\bold V(f))/\bold D(\Sigma)$ 
is affine, hence $H^i (\bold P\setminus X_f)=0$ for $i>d$.
We will prove by induction on $s$ that 
$H^i (\bold P\setminus(X_{f_1}\cap\ldots\cap X_{f_s}))=0$
for $i>d+s-1$.
Consider the Mayer-Vietoris sequence                      
$$\cdots\rightarrow H^i(U\cap V)\rightarrow 
H^{i+1}(U\cup V)\rightarrow H^{i+1}(U)\oplus H^{i+1}(V)\rightarrow
H^{i+1}(U\cap V)\rightarrow\cdots$$
with $U=\bold P\setminus(X_{f_1}\cap\ldots\cap X_{f_{s-1}})$,
 $V=\bold P\setminus X_{f_s}$.
Notice that $U\cup V=\bold P\setminus(X_{f_1}\cap\ldots\cap X_{f_s})$ and 
$U\cap V=\cup_{i=1}^{s-1}\bold P\setminus(X_{f_i}\cup
 X_{f_s})=\bold P\setminus(X_{f_{1}\cdot f_{s}}\cap\ldots\cap 
X_{f_{s-1}\cdot f_{s}})$.
So, using the induction and the above sequence, we obtain that 
$H^i(\bold P\setminus X)=0$ for $i>d+s-1$.
As a consequence of this, $X$ is nonempty unless $s>d$
because the dimension $h^{2d}(\bold P)=1$.
Since $\bold P\setminus X$ is a V-manifold, 
Poincar\'e duality implies that $H^{i}_{c}(\bold P\setminus X)=0$ for 
$i\leq d-s$.
Now the desired result follows from the long exact sequence of the 
cohomology with compact supports
($X$ and $\bold P$ are compact):
$$\cdots\rightarrow H^i_c(\bold P\setminus X)\rightarrow H^i_c(\bold P)\rightarrow
H^i_c(X)\rightarrow H_c^{i+1}(\bold P\setminus X)\rightarrow H_c^{i+1}(\bold P)\rightarrow\cdots.$$

If $X$ is an intersection of ample hypersurfaces defined by 
$f_1,\dots,f_s$,
then  Lemma 9.15 [BC] gives us that $f_1,\dots,f_s$ belong to $B(\Sigma)$.
\qed

{\bf Corollary 1.5.} {\it A quasi-smooth intersection 
$X=X_{f_1}\cap\ldots\cap X_{f_s}$,
defined by $f_1,\dots,f_s\in B(\Sigma)$,
has pure dimension $d-s$.}

Since the dimension of $H^0(X,\bold C)$ is 
the number of connected components of $X$, we obtain
another important  result.

{\bf Corollary 1.6.}
{\it An intersection $X_{f_1}\cap\ldots\cap X_{f_s}$,
defined by $f_1,\dots,f_s\in B(\Sigma)$,
 in a complete simplicial toric variety $\bold P_\Sigma$
is connected provided $s<\dim \bold P_\Sigma$.}

{\bf Remark 1.7.} If the polynomials                
 $f_1,\dots,f_s$ have ample degrees, then this corollary follows 
from a more general statement in [FL1] 
(see also [FL2] and [FH] for connectedness theorems).
              
{\bf 2. ``Cayley trick".}
We will explore a Cayley trick to reduce studying of the
cohomology of quasi-smooth intersections to results already known 
for hypersurfaces.                                            
                                                                               
Let $L_1,\dots,L_s$ be line bundles  on a complete $d$-dimensional
toric variety $\bold P=\bold P_\Sigma$, and let 
$\pi:\bold P(E)\rightarrow\bold P$ be the projective space bundle
associated to the vector bundle $E=L_1\oplus\cdots\oplus L_s$.
Then the $\Bbb P^{s-1}$-bundle $\bold P(E)$ is a toric variety.
The fan corresponding to it can be described as follows [O, p.~58].
Suppose that support functions $h_1,\dots,h_s$ give rise 
to the isomorphism classes of line bundles 
$[L_1],\dots,[L_s]\in\text{Pic}(\bold P)$, 
respectively.
Introduce  a $\bold Z$-module $N^\prime$ with a $\bold Z$-basis
$\{n_2,\dots,n_s\}$ and let $\tilde{N}:=N\oplus N^\prime$ and 
$n_1:=-n_2-\cdots-n_s$. Denote by $\tilde{\sigma}$ the image
of each $\sigma\in\Sigma$ under the $\bold R$-linear map 
$N_{\bold R}\rightarrow\tilde{N}_{\bold R}$
which sends $y\in N_{\bold R}$ to $y-\sum^{s}_{j=1}h_j(y)n_j$.
On the other hand, let $\sigma_i^\prime$ be the
cone in $N^\prime_{\bold R}$ generated by $n_1,\dots,n_i,n_{i+1},\dots,n_s$
and let $\Sigma^{\prime}$ be the fan in $N^\prime_{\bold R}$ 
consisting of the faces of $\sigma_1^\prime,\dots,\sigma_s^\prime$.
Then $\bold P(E)$ corresponds to the fan
$\tilde{\Sigma}:=\{\tilde\sigma+\sigma^{\prime}:\,\sigma\in\Sigma, 
\sigma^\prime\in\Sigma^\prime\}$.
From this  description it is easy to see that if
$\Sigma$ is a complete simplicial fan then
$\bold P(L_1\oplus\cdots\oplus L_s)$ is a complete 
simplicial toric variety.
We see that the integral generators of the 1-dimensional cones in
$\tilde{\Sigma}$ are given by 
$$\matrix
\tilde{e}_i =e_i-\sum\limits_{1\leq j\leq s}h_j(e_i)n_j,\quad 
i=1,\dots n,\\
\tilde{n}_1 =-n_2-\cdots-n_s,\hfill\\
\tilde{n}_j =n_j,\quad j=2,\dots,s,\hfill
\endmatrix$$
where $e_1,\dots,e_n$ are the integral generators of the 1-dimensional
cones in $\Sigma$.

The homogeneous coordinate ring of $\bold P(E)$ is the polynomial ring 
$$R=\bold C[x_1,\dots,x_n,y_1,\dots,y_s],$$
where $x_i$ corresponds to $\tilde{e}_i$ and $y_j$ corresponds to
$\tilde{n}_j$.
This ring has  a grading by the Chow group $A_{d+s-2}(\bold P(E))$.
Since $\bold P$ is a normal variety, there is an embedding of the Picard group
$\text{Pic}(\bold P)\hookrightarrow A_{d-1}(\bold P)$.
We want to show that if some polynomials 
$f_j\in S(\Sigma)=\bold C[x_1,\dots,x_n]$ have the 
property $\text{deg}(f_j)=[L_j]\in\text{Pic}(\bold P)$,
then the polynomials $y_{j}f_{j}$ all have the same degree in $R$.
This will allow us to consider a hypersurface defined by the homogeneous
polynomial $F=\sum_{j=1}^{s}y_{j}f_{j}$.
           
{\bf Lemma 2.1.}
{\it
Let $f_1,\dots,f_s\in S(\Sigma)$ be $\bold D$-homogeneous polynomials, 
such that $\deg(f_j)=[L_j]$ for some line bundles $L_1,\dots,L_s$.
Then  $F=\sum^{s}_{j=1}y_j f_j$ is
homogeneous in $R$ and its degree is the isomorphism class $[O_{\bold P(E)}(1)]$
of the canonical line bundle on $\bold P(E)=\bold P(L_1\oplus\cdots\oplus L_s)$.}

{\it Proof.} To prove that $F$ is a homogeneous polynomial we will repeat
the arguments in the proof of  Lemma 3.5 in [CCD].
Let $D_1,\dots,D_n$ be the torus-invariant divisors on $\bold P=\bold P_{\Sigma}$  
corresponding to the 1-dimensional cones of the fan $\Sigma$.
Then the pullback $\pi^* D_i$ is the torus-invariant divisor of $\bold P(E)$
corresponding to the cone generated by $\tilde{e}_i$. Also 
denote by $D_{j}^{\prime}$ the torus-invariant divisor corresponding to
$\tilde{n}_j$.
Let $\tilde{M}=M\oplus M^{\prime}$ be the lattice dual to 
$\tilde{N}=N\oplus N^\prime$ with $M^{\prime}=\text{Hom}(N^{\prime},\bold Z)$
having $\{n_2^*,\dots,n_s^*\}$ as a basis dual to $\{n_2,\dots,n_s\}$.
The divisor corresponding to the character $\chi^{n^*_j}$ is 
$$
\text{div}(\chi^{n_j^*})=\sum\limits^{n}_{i=1}\langle n^*_j,\tilde{e}_i\rangle
\pi^*D_i+\sum\limits^s_{k=1}\langle n^*_j,\tilde{n}_k\rangle D_k^{\prime}
=\sum\limits^{n}_{i=1}(h_1(e_i)-h_j(e_i))\pi^*D_i-D_1^{\prime}+D_j^{\prime}.$$
Therefore, $[D_j^{\prime}]+[\pi^*L_j]$ 
all have the same degree in the Chow
group $A_{d+s-2}(\bold P(E))$, and, consequently, $F$ is a homogeneous
polynomial. 

Now consider the following exact sequence [M]:
$$0\rightarrow O_{\bold P(E)}\rightarrow\pi^*E^*\otimes 
O_{\bold P(E)}(1)\rightarrow
\text{T}_{\bold P(E)}\rightarrow\pi^*\text{T}_{\bold P}\rightarrow 0,$$
where $\text{T}_X$ denotes the tangent bundle, $E^*$ is the dual bundle.
From here we can compute the Chern class 
$$c_1(\text{T}_{\bold P(E)})\!=\!
c_1(\pi^*\text{T}_{\bold P})+c_1(\pi^*E^*\otimes O_{\bold P(E)}(1))\!=\!
\pi^*c_1(\text{T}_{\bold P})-\pi^*c_1(E)+s\cdot c_1(O_{\bold P(E)}(1)).$$
Hence, $s\cdot c_1(O_{\bold P(E)}(1))=\pi^*c_1(L_1)+\cdots+\pi^*c_1(L_s)+
c_1(\text{T}_{\bold P(E)})-\pi^*c_1(\text{T}_{\bold P})$.
On the other hand, from the generalized Euler exact sequence [BC,~\S 12]
we get
$$0\rightarrow O_{\bold P}^{n-d}\rightarrow\oplus_{i=1}^{n}O_{\bold P}(D_i)\rightarrow
\text{T}_{\bold P}\rightarrow0.$$
This implies that $c_1(\text{T}_{\bold P})=[D_1]+\cdots+[D_n]$.
Similarly we have 
$c_1(\text{T}_{\bold P(E)})=[\pi^*D_1]+\cdots+[\pi^*D_n]+[D_1^\prime]+\cdots+[D_s^\prime]$.
Under the  identification 
$\text{Pic}(\bold P(E))\hookrightarrow A_{d+s-2}(\bold P(E))$
the first Chern class of a line bundle on $\bold P(E)$
is exactly its isomorphism class in the Picard group
$\text{Pic}(\bold P(E))$.
Therefore 
$$s\cdot[O_{\bold P(E)}(1)]=[\pi^*L_1]+\cdots+[\pi^*L_s]+[D_1^\prime]+
\cdots+[D_s^\prime]=s\cdot([\pi^*L_2]+[D_2^\prime]).$$
It can be easily checked that 
$D_2^\prime$ is a Cartier divisor on $\bold P(E)$.
Hence all classes $[O_{\bold P(E)}(1)]$, $[\pi^*L_2]$ and $[D_2^\prime]$
lie in the Picard group
$\text{Pic}(\bold P(E))$.
But this group is free abelian, because $\bold P(E)$ is complete.
So the above equality is divisible by $s$:
$[O_{\bold P(E)}(1)]=[\pi^*L_2]+[D_2^\prime]=\text{deg}(F)$.
\qed

From now on  we assume that 
$\bold P=\bold P_{\Sigma}$ is a complete simplicial toric variety
and that $\deg(f_j)\in\text{Pic}(\bold P)$, $j=1,\dots,s$.
Denote by $Y$ the hypersurface in $\bold P(E)$ defined by 
$F=\sum^{s}_{j=1}y_j f_j$.
           
{\bf Lemma 2.2.}
 {\it $X=X_{f_1}\cap\ldots\cap X_{f_s}$ is a quasi-smooth intersection
iff the hypersurface $Y$ is quasi-smooth.}

{\it Proof.}
 $X=X_{f_1}\cap\ldots\cap X_{f_s}$ is a quasi-smooth intersection
means that whenever $x\in\bold V(f_1,\ldots,f_s)\setminus Z(\Sigma)$,
the  $\text{rank}\bigl(\frac{\partial f_j}{\partial x_i}(x)\bigr)_{i,j}=s$.
And $Y$ is quasi-smooth iff $z=(x,y)\in\bold V(F)\setminus Z(\tilde\Sigma)$
implies that one of the partial derivatives 
$\frac{\partial F}{\partial y_j}(z)=f_j(x)$, $j=1,\dots,s$,
$\frac{\partial F}{\partial x_i}(z)=\sum^s_{j=1}y_j\frac{\partial f_j}
{\partial x_i}(x)$, $i=1,\dots,n$,
is nonzero.
 
So let $(x,y)\in\bold V(F)\setminus Z(\tilde\Sigma)$, then
there is  a cone $\tilde\sigma+\sigma^\prime\in\tilde\Sigma$
with $\sigma\in\Sigma$, $\sigma^\prime\in\Sigma^\prime$,
such that $\prod_{\tilde e_i\notin\tilde\sigma}x_i
\prod_{\tilde n_j\notin\sigma^\prime}y_j\ne0$
where $x_i, y_j$ are the coordinates of $(x,y)$.
If $f_1(x)=\cdots=f_s(x)=0$, then $x\in\bold V(f_1,\dots,f_s)\setminus 
Z(\Sigma)$ because $\prod_{e_i\in\sigma}x_i\ne0$.
And if $X=X_{f_1}\cap\ldots\cap X_{f_s}$ is 
a quasi-smooth intersection, one of the partial derivatives 
$\frac{\partial F}{\partial x_i}(z)=\sum^s_{j=1}y_j\frac{\partial f_j}
{\partial x_i}(x)$, $i=1,\dots,n$,
is nonzero.

Conversely, suppose $Y$ is quasi-smooth.
Pick any $x\in\bold V(f_1,\dots,f_s)\setminus Z(\Sigma)$,
then $(x,y)\in\bold V(F)\setminus Z(\tilde\Sigma)$
for each $y=(y_1,\dots,y_s)\ne0$.
Therefore  
$\sum^s_{j=1}y_j\frac{\partial f_j}
{\partial x_i}(x)\ne0$
for some $i$, which means the
$\text{rank}\bigl(\frac{\partial f_j}{\partial x_i}(x)\bigr)_{i,j}$
is maximal.\qed

{\bf 3. Cohomology of quasi-smooth intersections.}                  
Since a quasi-smooth intersection is a compact V-manifold (Proposition 1.3),
 the cohomology on it has a pure Hodge structure.                              
Using Proposition 1.4 and the Poincar\'e duality, we can compute
the cohomology of a quasi-smooth intersection except for the cohomology
in the middle dimension $d-s$.
So we introduce the following definition.

{\bf Definition 3.1.}
The {\it variable cohomology group} $H^{d-s}_{var}(X)$ is
$\text{coker}(H^{d-s}(\bold P)\overset{i^*}\to{\rightarrow}H^{d-s}(X)).$

The variable cohomology group also has a pure Hodge structure.

{\bf  Proposition 3.2.}
{\it Let $X=X_{f_1}\cap\ldots\cap X_{f_s}$ be a quasi-smooth intersection
of ample hypersurfaces.
Then there is an exact sequence of mixed Hodge structures
$$0\rightarrow H^{d-s-1}(\bold P)\overset{\cup[X]}\to{\rightarrow}
H^{d+s-1}(\bold P)\rightarrow H^{d+s-1}(\bold P\setminus X)\rightarrow
H^{d-s}_{var}(X){\rightarrow}0,$$
where $[X]\in H^{2s}(\bold P)$ is the cohomology class of $X$.}
          
{\it Proof.}  
Consider the Gysin exact sequence:                        
$$\cdots\rightarrow H^{i-2s}(X)\overset{i_!}\to{\rightarrow}
H^i(\bold P)\rightarrow H^i(\bold P\setminus X)\rightarrow
H^{i-2s+1}(X)\overset{i_!}\to{\rightarrow}H^{i+1}(\bold P)\rightarrow\cdots.\eqno{\text{(1)}}$$
Since $i^*$ is Poincar\'e dual to the Gysin map $i_!$,
it follows that $H_{var}^{d-s}(X)$ is isomorphic to the kernel
of $i_!:\,H^{d-s}(X)\rightarrow H^{d+s}(P)$.
So we get an exact sequence 
$$H^{d-s-1}(X)\overset{i_!}\to{\rightarrow}
H^{d+s-1}(\bold P)\rightarrow H^{d+s-1}(\bold P\setminus X)\rightarrow
H^{d-s}_{var}(X){\rightarrow}0.$$
Now we use a commutative diagram
$$\matrix
H^{d-s-1}(X) &\!\!\! \overset{i_!}\to{\rightarrow} H^{d+s-1}(\bold P)\hfill \\
{\scriptstyle i^*}\uparrow  &\!\!\! \nearrow_{\scriptstyle\cup[X]}\hfill  \\
H^{d-s-1}(\bold P). &                    
\endmatrix$$

By Proposition 1.4 $i^*$ is an isomorphism in this diagram, so it 
suffices to prove that the Gysin map $i_!$ 
is injective in the above diagram.
                     
{\bf Lemma 3.3.}
{\it If $X=X_{f_1}\cap\ldots\cap X_{f_s}$ is a quasi-smooth intersection
of ample hypersurfaces, then the Gysin map
$H^{d-s-1}(X)\overset{i_!}\to{\rightarrow}H^{d+s-1}(\bold P)$
is injective.}

{\it Proof.}
Since the odd dimensional cohomology of a complete 
simplicial toric variety vanishes [F1, pp.~92-94] and
$i^*:H^{d-s-1}(\bold P)\rightarrow H^{d-s-1}(X)$
is an isomorphism by Proposition 1.4, it follows
that $H^{d-s-1}(X)=H^{d-s-1}(\bold P)=H^{d+s-1}(\bold P)=0$
when $d+s-1$ is odd.
So by the Gysin exact sequence (1) it is enough to show that
$H^{d+s-2}(\bold P\setminus X)=0$ 
when $d+s-2$ is odd.
To prove this we use the Cayley trick again. Let $Y$ be the hypersurface
defined by $F=\sum^{s}_{j=1}y_j f_j$.
Then the natural map $\bold P(E)\setminus Y\rightarrow\bold P\setminus X$,
induced by the projection $\pi:\bold P(E)\rightarrow\bold P$,
is a $\bold C^{s-1}$ bundle in the Zariski topology.
Notice that $\bold P\setminus X$ is simply connected,
because $\bold P$ is simply connected [F1, p.~56]
and $X$ has codimension at least 2 in $\bold P$.
Hence, the Leray-Serre spectral sequence implies that                 
$H^i(\bold P(E)\setminus Y)=H^i(\bold P\setminus X)$
for $i\geq0$.
We have that $H^{d+s-2}(\bold P(E))=0$ for $d+s-2$ odd
and $Y$ is quasi-smooth by Lemma 2.2.
So from the Gysin exact sequence
$$H^{d+s-2}(\bold P(E))\rightarrow H^{d+s-2}(\bold P(E)\setminus Y)\rightarrow
H^{d+s-3}(Y)\overset{j_!}\to{\rightarrow}H^{d+s-1}(\bold P(E)).$$
(here the Gysin map $j_!$ is induced by the inclusion
$j:\,Y\hookrightarrow\bold P(E)$)
it follows that we need to show injectivity of
$j_!:H^{d+s-3}(Y)\rightarrow H^{d+s-1}(\bold P(E))$.
Consider the commutative diagram
$$\matrix
H^{d+s-3}(Y) &\!\!\! \overset{j_!}\to{\rightarrow} H^{d+s-1}(\bold P(E))\hfill \\
{\scriptstyle j^*}\uparrow  &\!\!\! \nearrow_{\scriptstyle\cup[Y]}\hfill  \\
H^{d+s-3}(\bold P(E)). &                    
\endmatrix$$
where $[Y]\in H^{2}(\bold P(E))$ is the cohomology class of $Y$.
The canonical line bundle $O_{\bold P(E)}(1)$ is ample [H,~III,\S 1],
whence by Lemma 2.1, $Y$ is ample.
So by Proposition 10.8 [BC]
$j^*:H^{d+s-3}(\bold P(E))\rightarrow H^{d+s-3}(Y)$
is an isomorphism and by Hard Lefschetz
$\cup[Y]:H^{d+s-3}(\bold P(E))\rightarrow H^{d+s-1}(\bold P(E))$
is injective.
Thus, from the above diagram the lemma follows.\qed

{\bf Definition 3.4.}
For a nonzero polynomial $F\in R=\bold C[x_1,\dots,x_n, y_1,\dots,y_s]$
the {\it Jacobian ring} $R(F)$ denotes the quotient of $R$ by the ideal
generated by the partial derivatives
$\frac{\partial F}{\partial y_j}$, $j=1,\dots,s$,
$\frac{\partial F}{\partial x_i}$, $i=1,\dots,n$.
            
{\bf Remark 3.5.}
If $F=y_1f_1+\cdots+y_sf_s$ is as in Lemma 2.1 with $f_j\in S_{\alpha_j}$,
 then $R(F)$
carries a natural grading by the Chow group $A_{d+s-2}(\bold P(E))$.
Moreover, there are canonical isomorphisms
$A_{d+s-2}(\bold P(E))\cong A_{d-1}(\bold P)\oplus A_d(\bold P)
\cong A_{d-1}(\bold P)\oplus \bold Z$
[F2].
With respect to this bigrading of the Chow group 
$A_{d+s-2}(\bold P(E))$ we have that $\deg(F)=(0,1)$,
$\deg(f_j)=(\alpha_j,0)$,  $\deg(y_j)=(-\alpha_j,1)$,
which is very similar to the case when $\bold P$ is a projective space.

We now can state the main result.
             
{\bf Theorem 3.6.}
{\it Let $\bold P$ be a $d$-dimensional complete simplicial toric variety,
and let $X\subset\bold P$ be a quasi-smooth intersection of 
ample hypersurfaces 
defined by $f_j\in S_{\alpha_j}$, $j=1,\dots,s$. If $F=y_1f_1+\cdots+y_sf_s$,
then for $p\ne\frac{d+s-1}{2}$, we have a canonical isomorphism
$$R(F)_{(d+s-p)\beta-\beta_0}\cong H_{var}^{p-s,d-p}(X)$$
where $\beta_0=\deg(x_1\cdots x_n\cdot y_1\cdots y_s)$,
$\beta=\deg(F)=\deg(f_j)+\deg(y_j)$.  
In the case $p=\frac{d+s-1}{2}$ there is  an exact sequence
$$0\rightarrow H^{d-s-1}(\bold P)\overset{\cup[X]}\to{\rightarrow}
H^{d+s-1}(\bold P)\rightarrow R(F)_{\frac{d+s+1}{2}\beta-\beta_0}\rightarrow
H^{\frac{d-s-1}{2},\frac{d-s+1}{2}}_{var}(X){\rightarrow}0.$$}

{\it Proof.}
Since $H^i(\bold P)$ vanishes for $i$ odd and has a pure Hodge structure
of type $(p,p)$ for $i$ even,                                
from Proposition 3.2 we get
$\text{Gr}_F^p H^{d+s-1}(\bold P\setminus X)\cong H_{var}^{p-s,d-p}(X)$
if $p\ne\frac{d+s-1}{2}$, and in case $p=\frac{d+s-1}{2}$ the following
sequence
$$0\!\rightarrow\! H^{d-s-1}(\bold P)\overset{\cup[X]}\to{\rightarrow}
H^{d+s-1}(\bold P)\!\rightarrow\!\text{Gr}_F^{\frac{d+s-1}{2}}H^{d+s-1}
(\bold P\setminus X)\rightarrow
H^{\frac{d-s-1}{2},\frac{d-s+1}{2}}_{var}(X){\rightarrow}0\hfill$$
is exact.

Now use the isomorphism of mixed Hodge structures 
$H^i(\bold P\setminus X)\cong H^i(\bold P(E)\setminus Y)$
and by the Theorem 10.6 [BC] the desired result follows.\qed
     
{\bf 4. Cohomology of nondegenerate intersections.}
In this section we consider a special case of quasi-smooth
intersections. 
                                       
{\bf Definition 4.1.}                           
A closed subset $X=X_{f_1}\cap\ldots\cap X_{f_s}$,
defined by $\bold D$-homogeneous polynomials
$f_1,\dots,f_s$, is called a {\it nondegenerate  
intersection} if $X_{f_{j_1}}\cap\ldots\cap X_{f_{j_k}}\cap\bold T_{\tau}$
is a smooth subvariety of codimension $k$ in $\bold T_\tau$
for any $\{j_1,\dots,j_k\}\subset\{1,\dots,s\}$ and $\tau\in\Sigma$.
(Here $\bold T_\tau$ denotes the torus in $\bold P_\Sigma$ associated
with a cone $\tau\in\Sigma$).                        
                            
We will show  how to define a nondegenerate intersection in terms of
the polynomials $f_1,\dots,f_s$.
For $\sigma\in\Sigma$, let $U_\sigma=\{x\in\bold A^n:\,\hat x_\sigma 
\ne0\}$. We know that $\bold P_\Sigma$ has an affine toric open cover
by $\bold A_\sigma=U_\sigma/\bold D(\Sigma)$, $\sigma\in\Sigma$ [BC].
Also $\bold T_\tau=(U_\tau\setminus\cup_{\gamma\prec\tau}U_\gamma)/
\bold D(\Sigma)$.
Notice that $U_\tau\setminus\cup_{\gamma\prec\tau}U_\gamma=\{x\in\bold A^n:
\,\hat x_\tau\ne0, x_i=0\text{ if } \rho_i\subset\tau\}$ is a torus.
So each $\bold T_\tau$ is a quotient of a torus by a $D$-subgroup, 
because $\bold D$ is diagonalizable  [BC]. 

{\bf Lemma 4.2.}
{\it Let $T=(\bold C^*)^n/G$ be the quotient of a torus by a $D$-subgroup
$G$. Suppose that $X\subset(\bold C^*)^n$ is an invariant subvariety
with respect to the action of $G$. Then the geometric quotient
$X/G$ is smooth iff $X$ is smooth.}

{\it Proof.}
By the structure theorem of a $D$-group [Hu,~\S 16.2]
we can assume that $(\bold C^*)^n=G^\circ\times(\bold C^*)^k$,
where $G^\circ\cong(\bold C^*)^{n-k}$ is the identity component
of $G$, and $G=G^\circ\times H$ for some finite subgroup $H$ in
$(\bold C^*)^k$. Now it suffices to show  the Lemma if $G$ is a torus
or a finite group.
If $G=G^\circ$ then $X=(\bold C^*)^{n-k}\times p(X)$,
where by $p(X)$ we mean the projection of $X$ onto $(\bold C^*)^k$.
Notice that $p(X)\cong X/G$, hence $X$ is smooth iff $X/G$ is smooth.
In the case $G=H$ is a finite group it can be easily checked that
$X\rightarrow X/G$ is an unramified cover [Sh,~p. 346]. 
So $X$ and $X/G$ are smooth simultaneously.\qed

From this Lemma it follows that $X=X_{f_1}\cap\ldots\cap X_{f_s}$
is a nondegenerate intersection iff $\bold V(f_{j_1},\dots,f_{j_k})\cap V_\tau$
is a smooth subvariety of codimension $k$ in the torus 
$V_\tau=\{x\in\bold A^n:\,\hat x_\tau\ne0,
x_i=0\text{ if }\rho_i\subset \tau\}$.

As in section 2 we can consider the hypersurface 
$Y\subset\bold P(E)$ defined by $F=\sum^s_{j=1}y_j f_j$.

{\bf Lemma 4.3.} 
 {\it $X=X_{f_1}\cap\ldots\cap X_{f_s}$ is a nondegenerate intersection
iff $Y$ is a nondegenerate hypersurface.}

{\it Proof.}
As  shown above, $X=X_{f_1}\cap\ldots\cap X_{f_s}$ is a nondegenerate
intersection if
the  $\text{rank}\bigl(\frac{\partial f_j}{\partial x_i}(x)\bigr)_
{i\in\{i:\,e_i\notin\tau\}}^{j\in\{j_1,\dots,j_k\}}=k$
 for all $x\in\bold V(f_{j_1},\dots,f_{j_k})\cap V_\tau$,
$\tau\in\Sigma$ and $\{j_1,\dots,j_k\}\subset\{1,\dots,s\}$.
Similarly $Y$ is nondegenerate iff 
$z=(x,y)\in\bold V(F)\cap V_{\tilde{\tau}+\tau^\prime}$,
$\tilde{\tau}+\tau^\prime\in\tilde\Sigma$ 
with $\tau\in\Sigma$, $\tau^\prime\in\Sigma^\prime$ 
(recall the definition of $\bold P(E)$ 
associated with $\tilde\Sigma$ in the section 2) implies that
 one of the partial derivatives 
$\frac{\partial F}{\partial y_j}(z)=f_j(x)$, $j\in\{j:\,\tilde n_j\notin
\tau^\prime\}$, 
$\frac{\partial F}{\partial x_i}(z)=\sum^s_{j=1}y_j\frac{\partial f_j}
{\partial x_i}(x)$, $i\in\{i:\,\tilde{e}_i\notin\tilde{\tau}\}$,
is nonzero.
 
Let $(x,y)\in\bold V(F)\cap V_{\tilde\tau+\tau^\prime}$,
where $\tilde{\tau}+\tau^\prime\in\tilde\Sigma$ 
with $\tau\in\Sigma$, $\tau^\prime\in\Sigma^\prime$.
Then 
$\prod_{\tilde e_i\notin\tilde\tau}x_i
\prod_{\tilde n_j\notin\tau^\prime}y_j\ne0$
and $x_i=0$ if $\tilde e_i\in\tilde\tau$,
 $y_j=0$ if $\tilde n_j\in\tau^\prime$.
If $f_j(x)=0$ for all $j\in\{j:\,\tilde n_j\notin\tau^\prime\}$,
 then $x\in\bold V(f_{j_1},\dots,f_{j_k})\cap V_\tau$
where $\{j_1,\dots,j_k\}=\{j:\,\tilde n_j\notin\tau^\prime\}$.
So if $X=X_{f_1}\cap\ldots\cap X_{f_s}$ is 
a nondegenerate intersection, one of the partial derivatives 
$\frac{\partial F}{\partial x_i}(z)=\sum^s_{j=1}y_j\frac{\partial f_j}
{\partial x_i}(x)$, $i\in\{i:\tilde e_i\notin\tilde\tau\}$,
is nonzero.

Conversely, suppose $Y$ is  nondegenerate.
Take any $x\in\bold V(f_{j_1},\dots,f_{j_k})\cap V_\tau$ with
$\tau\in\Sigma$, $\{j_1,\dots,j_k\}\subset\{1,\dots,s\}$.
Then $(x,y)\in\bold V(F)\cap V_{\tilde\tau+\tau^\prime}$
for each $y\in V_{\tau^\prime}=\{y\in\bold A^s:\,y_j\ne0\text{ if }
\tilde n_j\notin\tau^\prime, y_j=0\text{ if }\tilde n_j\in\tau^\prime\}$
where $\tau^\prime$ is the cone generated by the complement of 
$\{\tilde n_{j_1},\dots,\tilde n_{j_k}\}$ in the set 
$\{\tilde n_1,\dots,\tilde n_s\}$.
Therefore  
$\sum^s_{j=1}y_j\frac{\partial f_j}
{\partial x_i}(x)\ne0$
for some $i$, which means the
$\text{rank}\bigl(\frac{\partial f_j}{\partial x_i}(x)\bigr)_
{i\in\{i:\,e_i\notin\tau\}}^{j\in\{j_1,\dots,j_k\}}=k$.\qed

Since a nondegenerate hypersurface is quasi-smooth [BC],
 Lemma 2.2 shows that a nondegenerate intersection
is quasi-smooth.

{\bf Definition 4.4.} [BC]
Given a polynomial $f\in S=\bold C[x_1,\dots,x_n]$,
we get the ideal quotient 
 $J_1(f)=\langle x_1\partial f/\partial x_1,\dots,x_n
\partial f/\partial x_n\rangle:x_1\cdots x_n$ (see [CLO,~p.~193]) 
and the ring $R_1(f)=S/J_1(f)$.

{\bf Remark 4.5.}
If $F=\sum^s_{j=1}y_j f_j\in R$ is as in Lemma 2.1, then 
$R_1(F)=R/J_1(F)$ has a natural grading  by the Chow group 
$A_{d+s-2}(\bold P(E))\cong A_{d-1}(\bold P)\oplus\bold Z$.

{\bf Theorem 4.6.}
{\it Let $X=X_{f_1}\cap\ldots\cap X_{f_s}$ be a nondegenerate intersection
of ample hypersurfaces given by $f_j\in S_{\alpha_j}$, $j=1,\dots,s$.
If $F=\sum^s_{j=1}y_j f_j\in R$, then there is a canonical
isomorphism
$$H^{p-s,d-p}_{var}(X)=R_1(F)_{(d+s-p)\beta-\beta_0},$$
where $\beta_0=\deg(x_1\cdots x_n\cdot y_1\cdots y_s)$,
$\beta=\deg(F)$.}  
                              
{\it Proof.}
First we will show that there is an isomorphism of Hodge structures
$H^{d-s}_{var}(X)(1-s)\cong H^{d+s-2}_{var}(Y)$.
Let $\varphi:Y\rightarrow\bold P$ be the composition of the inclusion
$j:\,Y\hookrightarrow\bold P(E)$ and the projection 
$\pi:\bold P(E)\rightarrow\bold P$.
As in [Te], consider the following morphism of the Leray spectral sequences
$$\matrix
E_2^{p,q}=& H^p(\bold P, R^q\pi_*\bold C) & \Rightarrow &H^{p+q}(\bold P(E)) \\
 &\downarrow&&\downarrow\\
^\prime E_2^{p,q}=& H^p(\bold P, R^q\varphi_*\bold C)&\Rightarrow & H^{p+q}(Y).
\endmatrix$$
Since 
$$\varphi^{-1}(X)=\left\{\matrix \Bbb P^{s-1}\text{ if }x\in X,
\\ \Bbb P^{s-2}\text{ if }x\notin X,
\endmatrix\right.$$
we have that (see [Go,~p.~202], [De])
$$R^q\varphi_*\bold C=\left\{\matrix \bold C_{\bold P}(-\frac{q}{2})\hfill &
\text{ if }q\text{ is even and }0\leq q<2s-2,\hfill\\
\bold C_X(1-s) & \text{ if }q=2s-2,\hfill\\
0 & \text{ otherwise}.\hfill
\endmatrix\right.$$
Also we have 
$$R^q\pi_*\bold C=\left\{\matrix \bold C_{\bold P}(-\frac{q}{2}) &
\text{ if }q\text{ is even and }0\leq q\leq 2s-2,\hfill\\
0 & \text{ otherwise}.\hfill
\endmatrix\right.$$

The first spectral sequence degenerates at $E_2$, 
because for $p$ or $q$ odd $E_r^{p,q}$ vanishes.
The second spectral sequence also degenerates at $E_2$:
$$h^{l-2s-2}(X)+\sum_{q=0}^{2s-4}h^{l-q}(\bold P)=\sum_{p+q=l}\dim
{}^\prime E_2^{p,q}\geq
\sum_{p+q=l}\dim {}^\prime E_{\infty}^{p,q}=h^l(Y).$$
To show the degeneracy of $^\prime E_2^{p,q}$
it suffices to show that the above inequality is an equality.
From Proposition 10.8 [BC] and Proposition 3.2 we get 
$$h^{d+s-2}(Y)=h^{d+s-2}(\bold P(E))+h^{d+s-1}(\bold P(E)\setminus Y)-h^{d+s-1}(\bold P(E))
+h^{d+s-3}(\bold P(E)),$$
$$h^{d-s}(X)=h^{d-s}(\bold P)+h^{d+s-1}(\bold P\setminus X)-h^{d+s-1}(\bold P)
+h^{d-s-1}(\bold P).$$
Hence, using the spectral sequence $E_2^{p,q}$, we can easily 
compute the Hodge numbers of $\bold P(E)$ and check that 
$h^{l-2s-2}(X)+\sum_{q=0}^{2s-4}h^{l-q}(\bold P)=h^l(Y)$
for $l=d+s-2$.
Using Proposition 1.4, we can similarly show the above equality
 for $l\ne d+s-2$ as well.
So the spectral sequence $^\prime E_2^{p,q}$ degenerates at $E_2$.
Since $E_2^{d+s-2-q,q}={}^\prime E_2^{d+s-2-q,q}$ for $q\ne 2s-2$ 
and, by Proposition 1.4,
$E_2^{d-s,2s-2}\hookrightarrow {}^\prime E_2^{d-s,2s-2}$,
we get an isomorphism of Hodge structures (for details see [Te]):
$$H^{d+s-2}_{var}(Y)\cong {}^\prime E^{d-s,2s-s}_2/E_2^{d-s,2s-2}
\cong H^{d-s}_{var}(X)(1-s).$$
Now we only need to apply Theorem 11.8 [BC]
to finish the proof.\qed

\enddocument